\documentclass[12pt]{article}
\usepackage{amsmath}
\usepackage{amsfonts}
\usepackage{amssymb}
\usepackage{latexsym}
\usepackage{graphicx}
\usepackage{enumerate}
\usepackage[shortlabels]{enumitem}
\usepackage{mathtools}
\usepackage{hyperref}
\usepackage{arydshln}

\usepackage{amsmath, amssymb,latexsym}
\usepackage{color}
\usepackage{xcolor}
\usepackage{fullpage}

\def\C{\mathbb{C}}

\def\la{\lambda}
\def\N{\mathcal{N}}
\def\T{\mathsf{*}}
\def\V{\mathcal{V}}
\def\U{\mathcal{U}}

\DeclareMathOperator{\diag}{diag}

\DeclareMathOperator{\rank}{rank}
\DeclareMathOperator{\im}{Im}

\newcommand{\tm}{\color{blue} \times}
\newcommand{\tmg}{\color{gray} \times}
\newcommand{\tmr}{{\color{red} \times}}
\newcommand{\otb}{\color{blue} \otimes}
\newcommand{\otr}{\color{red} \otimes}
\newcommand{\zg}{\color{gray} 0}
\newcommand{\zb}{\color{blue} 0}

\newcommand{\zr}{\color{red} 0}

\newtheorem{theorem}{Theorem}[section]
\newtheorem{proposition}[theorem]{Proposition}

\newtheorem{lemma}[theorem]{Lemma}
\newtheorem{definition}[theorem]{Definition}
\newtheorem{corollary}[theorem]{Corollary}
\newtheorem{remark}[theorem]{{\sc Remark}}
\newtheorem{example}[theorem]{Example}

\title{On computing root polynomials and minimal bases of matrix pencils}
\date{}
\author{
Vanni Noferini\thanks{Aalto University, Department of Mathematics and Systems Analysis, P.O. Box 11100, FI-00076, Aalto, Finland. Supported by an Academy of Finland grant (Suomen Akatemian p\"{a}\"{a}t\"{o}s 331240). Email: \texttt{vanni.noferini@aalto.fi}}
\and 
Paul Van Dooren\thanks{Universit\'e Catholique de Louvain, Department of Mathematical Engineering, Av. Lemaitre 4, B-1348 Louvain-la-Neuve, Belgium. Supported by an Aalto Science Institute Visitor Programme. Email: \texttt{paul.vandooren@uclouvain.be}}
}
\begin{document}
\maketitle
\begin{abstract}
 We revisit the notion of root polynomials, thoroughly studied in [F. Dopico and V. Noferini, Root polynomials and their role in the theory of matrix polynomials, Linear Algebra Appl. 584:37--78, 2020] for general polynomial matrices, and show how they can efficiently be computed in the case of a matrix pencil $\lambda E +A$. The method we propose makes extensive use of the staircase algorithm, which is known to compute the left and right minimal indices of the Kronecker structure of the pencil. In addition, we show here that the staircase algorithm, applied to the expansion $(\lambda-\lambda_0)E+(A-\la_0 E)$, constructs a block triangular pencil from which a minimal basis and a maximal set of root polynomials at the eigenvalue $\la_0$, can be computed in an efficient manner. 
\end{abstract}

\textbf{Keywords:} Root polynomial, maximal set, minimal basis, matrix pencil, staircase algorithm, Smith form, local Smith form

\textbf{MSC:} 15A03, 15A09, 15A18, 15A21, 15A22, 65F15

\section{Introduction}

Finding the eigenvalues of a polynomial matrix
\begin{equation*}
P(\la) = P_0 + P_1 \la + \ldots + P_d \la^d \; \in \C[\lambda]^{m\times n},
\end{equation*} 
and their partial multiplicities is a problem that occurs
naturally when one wants to describe the solution set of particular
matrix equations involving polynomial matrices \cite{GLR82,kailath,Kar94,kucera}. From the theoretical point of view, the problem is completely solved by the existence of a local Smith form at any point $\la_0 \in \C$ \cite{GLR82}. Computationally, the local Smith form can be determined by finding certain important polynomial vectors associated with $P(\la)$, namely, left and right \emph{minimal bases} \cite{For75, Mac21} and \emph{root polynomials} \cite{DopN21, GLR82, N11}. It has been long known that vectors in a minimal basis carry the information on the minimal indices \cite{For75}. That root polynomials have a similarly important role has been recently advocated: for example, so-called maximal sets of root polynomials encode all the information on partial multiplicities \cite{DopN21}, and they make it possible to properly define eigenvectors for singular polynomials \cite{DopN21}, which has proved useful for instance to carry out probabilistic studies of the condition number of eigenvalues \cite{LN20}.

This motivates a natural algorithmic question: how can one 
compute the vectors in a minimal basis and a maximal set of root polynomials of a polynomial matrix? This paper focuses on the case of $d=1$, i.e., matrix pencils, and addresses the question by providing a robust algorithm which builds on the staircase algorithm \cite{vd79}. There are at least two reasons to give special attention to pencils: first, the generalized eigenvalue problem is arguably the most common instance of polynomial eigenvalue problems; and second, even if the eigenvalue problem has higher degree to start with, it is common practice to linearize it as a first step towards its numerical solution. For many commonly used linearizations, it is known both how to recover minimal bases (a topic addressed by many papers in the literature, each focusing on different classes of linearizations: see for example \cite{DDM09,DDM10,DLPVD,NP16}) and how to recover maximal sets of root polynomials (a problem solved in \cite{DopN21} for many classes of linearizations) for the linearized polynomial matrix, starting from their counterparts for the linearizing pencils. Hence, in this precise sense an algorithm that solves the problem for pencils can be easily extended to an algorithm that computes root polynomials and minimal bases for a polynomial matrix of any degree.

The structure of the paper is as follows. In Section \ref{Sec:Background} we recall the necessary background and definitions
for zero directions, root polynomials and minimal bases of an arbitrary polynomial matrix. In Section \ref{Sec:Pencils} we consider the special case of a matrix pencil and show the link between zero directions, the rank profile of certain bidiagonal block Toeplitz matrices, and the construction of a so-called Wong sequence of subspaces. We also recall how the staircase algorithm for pencils of matrices constructs particular bases for such a Wong sequence. In Section \ref{Sec:Extraction} we show how a particular bidiagonalization procedure allows us to extract from this a maximal set of root polynomials, on one hand, and a minimal basis for the right null space, on the other hand. In Section \ref{Sec:Recurrences} we then develop simple recurrences that compute, for a given pencil, a minimal basis for the right null space, and a maximal complete set of $\la_0$-independent root polynomials.  In Section \ref{Sec:Numerics} we give numerical examples to illustrate our algorithm and we comment on the computational complexity. We end with a few concluding remarks in Section \ref{Sec:Conclusion}.

\section{Background and definitions}  \label{Sec:Background}

\subsection{Zeros, the local Smith form and null spaces}

 A finite zero, or eigenvalue, of $P(\la)$ is an element $\lambda_0 \in \C$ such that 
\[ \mathrm{rank}_\C P(\lambda_0) < \mathrm{rank}_{\C(\la)} P(\la).\]
The structure at a point $\lambda_0$ which is
a zero of a polynomial matrix $P(\lambda)$ is defined via
the local Smith form of the $m \times n$ polynomial matrix
$P(\lambda)$ at the point $\lambda_0\in \C$:
\begin{equation} \label{Smith}
M(\lambda)\cdot P(\lambda)\cdot N(\lambda) :=
\left[ \begin{array}{ccc|c} (\lambda-\lambda_0)^{\sigma_1} & & 0 & \\
 & \ddots & & \\ 0 & &  (\lambda-\lambda_0)^{\sigma_r}& \\ \hline
 & & & 0_{m-r,n-r} \end{array}\right], 
\end{equation}
and where $M(\lambda)$ and $N(\lambda)$ are polynomial and invertible at
$\lambda_0$, whereas $r=\mathrm{rank}_{\C(\la)} P(\la)$ is the normal rank of $P(\la)$. Furthermore, the integers $\sigma_i$ are known as the \emph{partial multiplicities}, or {\it
structural indices}, of $P(\lambda)$ at the zero
$\lambda_0$; they satisfy $0\leq \sigma_1 \leq \ldots \leq \sigma_r \leq dr$. The finite sequence
$\sigma_1, \ldots, \sigma_r$, or rather its subsequence listing its positive elements, is sometimes also called the Segr\'e characteristic 
at $\lambda_0$. The classical algorithm for the computation of the above decomposition is
based on the Euclidean algorithm and Gaussian elimination over the
ring of polynomials, which is in general numerically unreliable
\cite{Zuniga}.
For this reason it can be replaced by a technique, based on the
expansion around the point $\lambda_0$ \cite{vdv}, as explained in 
the next sections.
 
Other important sets of indices of a general $m\times n$ polynomial matrix 
$P(\la)$ are related to its right nullspace $\N_r$ and left nullspace $\N_\ell$,
which are rational vector spaces over the field $\C(\la)$ of rational functions in $\la$.
For this, we first need the following definition.
\begin{definition} 
The columns of a polynomial matrix $N(\la) \in \C[\lambda]^{n\times p}$ of normal rank $p$
is called a minimal polynomial basis if the sum of the degrees of its columns, 
called the order of the basis, is the minimal among all bases of 
span $N(\la)$. Its ordered column degrees are called the minimal indices of the basis.
\end{definition}
A minimal basis is not uniquely determined by the subspace it spans, but it was shown in \cite{For75} that the minimal indices are. If we define the right nullspace $\N_r(P)$ and the left nullspace 
$\N_\ell(P)$ of an $m\times n$  polynomial matrix $P(\la)$ of normal rank $r$ as 
the vector spaces of rational vectors $x(\la)$ and $y(\la)$ annihilated by $P(\la)$
$$  \N_r(P):=\{ x(\la) \; | \; P(\la)x(\la)=0\}, \qquad    \N_\ell(P) :=\{ y(\la) \;| \; y^\T(\la) P(\la)=0\}
$$
then the minimal indices of any minimal polynomial basis for these spaces, are called the right and left  
minimal indices of $P(\la)$. Their respective dimensions are $n-r$ and $m-r$ and the respective 
indices are denoted by
$$\{\epsilon_1,\ldots,\epsilon_{n-r}\}, \quad \{\eta_1,\ldots,\eta_{m-r}\}.
$$
It was also shown in \cite{For75} that, for any minimal basis $N(\la)$, the constant matrix 
$N(\la_0)$ has full column rank for all $\la_0\in\C$ and the highest column degree matrix of $N(\la)$ (or equivalently the columnwise reversal of $N(\la)$ evaluated at $\la=0$ \cite{NP}) also has full column rank.

\subsection{Zero directions, root polynomials and null vectors}

Let us assume that $\lambda_0$ is a zero of $P(\lambda)$ and let us express the latter polynomial matrix by its Taylor expansion around $\la_0$~:
\begin{equation*} \label{Toeplitz}
P(\lambda) := P_0 + (\lambda-\lambda_0)P_1 +
(\lambda-\lambda_0)^2P_2 + \ldots + (\lambda-\lambda_0)^dP_d .
\end{equation*}
Let us then define for $k\ge 0$ the following associated Toeplitz matrices and their ranks as
\begin{equation*} T_{\lambda_0,k} := \left[ \begin{array}{cccc}
P_0 & P_1 & \ldots & P_k \\
& P_0 & \ddots & \vdots \\
& & \ddots & P_1 \\
& & & P_0
\end{array}\right], \quad r_k:= \rank T_{\la_0,k}, 
\end{equation*}
where we implicitly have set the coefficients $P_i=:0$ for $i>d$.
Below, we will drop the suffix $\lambda_0$ when it is obvious from the context that we use an expansion about that point.
The rank increments $\rho_j := r_{j} - r_{j-1}$ were shown in \cite{vdv} to completely determine the partial multiplicities of $\la_0$, $\sigma_i (i=1,\ldots,r)$  and we can thus expect that the definition of root polynomials also should be related.

\medskip

Throughout this paper, we extend the usual notation of modular arithmetic from scalars to matrices by applying it elementwise. Namely, given a scalar polynomial $p(\lambda) \in \C[\la]$ and two polynomial matrices of the same size, say, $A(\la),B(\la) \in \C[\la]^{m \times n}$, then the notation $A(\la) \equiv B(\la) \mod p(\la)$ is shorthand to mean that there  exists a third polynomial matrix $C(\la) \in \C[\la]^{m \times n}$ such that $A(\la)-B(\la)=p(\la) C(\la)$. Therefore, for example, $\begin{bmatrix}
\la^2 & \lambda^2-3\la+2
\end{bmatrix} \equiv \begin{bmatrix}
\la & 0
\end{bmatrix} \mod (\la-1)$.

Let now $x(\lambda)\in \C[\lambda]^n$ and $y(\lambda)\in \C[\lambda]^m$ be polynomial vectors satisfying

\begin{eqnarray*}  x(\lambda)  \equiv x_0 + (\lambda-\lambda_0)x_1 +
\ldots + (\lambda-\lambda_0)^{k-1}x_{k-1} &\mod (\la-\la_0)^k, \quad x(\lambda_0)=x_0 \neq 0\\
  y(\lambda) \equiv y_0 + (\lambda-\lambda_0)y_1 +
\ldots + (\lambda-\lambda_0)^{k-1}y_{k-1} &\mod (\la-\la_0)^k,  \quad y(\lambda_0)=y_0 \neq 0
\end{eqnarray*}
then for $k>0$, we say that the vectors $x(\la)$ and $y(\la)$, respectively, are right and left {\it zero directions of order $k$} if
\begin{equation}\label{eq:sei} P(\lambda)x(\lambda) =  (\lambda-\lambda_0)^{k}v(\lambda), \quad v(\lambda) \in \C(\lambda)^m, \;
v(\lambda_0) \neq 0
\end{equation}
\begin{equation}\label{eq:sette}
 y^\T(\lambda)P(\lambda) =  (\lambda-\lambda_0)^{k}w(\lambda)^\T, \quad w(\lambda) \in \C(\lambda)^n, \;
w(\lambda_0) \neq 0.
\end{equation}
The fact that the respective vectors $x_0$ and $y_0$ are nonzero avoids trivial solutions \cite{EliK98,Kar94,VDooren93} obtained by multiplication 
with a positive power of $(\la-\la_0)$. Similarly, one can define right zero directions of order $\infty$ at $\la_0$ as polynomial vectors $x(\la)$ such that $x(\la_0)\neq0$ and $P(\la)x(\la)=0$; left zero directions of order $\infty$ are defined analogously.

Using \eqref{Smith}, expanding \eqref{eq:sei} and \eqref{eq:sette} in powers of $(\la-\la_0)$, and observing that the coefficient $(\la-\la_0)^j$ must be $0$ for all $j \geq k-1$, we obtain
\begin{equation*} 
T_{k-1}
\left[\begin{array}{c} x_{k-1} \\ \vdots \\ x_0 \end{array}\right]=0, \quad
\left[\begin{array}{ccc} y^\T_{0} & \ldots & y^\T_{k-1} \end{array}\right]
T_{k-1}=0.
\end{equation*}
So the problem of finding zero directions is apparently solved by computing the null spaces of these Toeplitz matrices \cite{VDooren93}. Unfortunately, though, knowledge of zero directions alone is not sufficient to extract the information about the minimal indices and the partial multiplicities at the point $\lambda_0$. However, there are two important subsets of zero directions of a polynomial matrix that allow us to do this, and thus deserve more attention.
For this we restrict ourselves to the right zero directions since the problem for the left zero directions is obtained by just
considering the conjugate transposed matrix $P(\la)^\T$.

For a singular polynomial matrix the set of zero directions also contains the polynomial vectors in the right nullspace
$\N_r(P)$. This follows easily from the fact that if $x(\la)$ is a column of a minimal basis matrix $N(\la)$, then $P(\la)x(\la)=0$ and $x(\lambda_0)\neq 0$ for any $\lambda_0\in \C$. It is thus a zero direction of order $\infty$ for any point 
$\la_0\in \C$. Moreover, if these zero directions are indeed the vectors of a minimal bases, their degrees provide all the information about the minimal indices.

Another special subset of zero directions is the set of root polynomials at a point $\la_0\in \C$. Root polynomials and their properties were studied in detail in \cite{DopN21}, where it was advocated that they deserve an important role in the theory of polynomial matrices; they had previously appeared as technical tools for proving other results \cite{GLR82, N11}. Below we give a definition which is clearly equivalent to that given in \cite{DopN21}, but rephrased in a way more convenient for us in this paper.
\begin{definition} Let the columns of $N(\la)$ be a right minimal basis of $P(\la)$; then
\begin{itemize} \item $r(\la)$ is a root polynomial of order $k$ if it is a zero direction of order $k$ and $\left[ N(\la_0) \; r(\la_0)\right]$ has full column rank
 \item $\{r_1(\la),...,r_s(\la)\}$ is a set of $\la_0$-independent root polynomials of orders $k_1,...,k_s$ if
they are zero directions of orders  $k_1,...,k_s$ and $\left[ N(\la_0) \; r_1(\la_0) \ldots r_s(\la_0) \right]$ has full column rank.
\item a $\la_0$-independent set is complete if there does not exist any larger $\la_0$-independent set
\item such a complete set is ordered if $k_1\ge \ldots \ge k_s > 0$
\item such a complete ordered set is maximal if there is no root polynomial $\tilde r(\la)$ of order $k>k_j$ at 
$\la_0$ such that $\left[ N(\la_0) \; r_1(\la_0) \ldots r_{j-1}(\la_0) \; \tilde r(\la_0) \right]$ has full column rank, for all $1\le j\le s$.
\end{itemize}
\end{definition}

The importance of maximal sets of root polynomials is given by the following result \cite[Theorem 4.1.3]{DopN21}:

\begin{theorem}\label{thm:frovan}
Let the nonzero partial multiplicities at $\lambda_0$ of $P(\la)$ be $0<\sigma_1 \leq \dots \leq \sigma_s$ and suppose that $r_1(\la),\dots,r_s(\la)$ are a complete set of root polynomials at $\la_0$ for $P(\la)$.  Then, the following are equivalent:
\begin{enumerate}
\item $r_1(\la),\dots,r_s(\la)$ are a maximal set of root polynomials at $\lambda_0$ for $P(\la)$;
\item the orders of such a set are precisely $\sigma_1,\dots,\sigma_s$;
\item the sum of the orders of such a set is precisely $\sum_{i=1}^s \sigma_i$.
\end{enumerate}
\end{theorem}

If the polynomial matrix $P(\la)$ is regular, then any zero direction is a root polynomial, and the definition of a maximal set can be applied to them too. However, in the singular case, generally zero directions do not provide the correct information on minimal indices and partial multiplicities. We illustrate this fact with the next simple example.

\begin{example}
The polynomial matrix
\[ P(\la) = \begin{bmatrix}
\la & \la\\
\la & \la
\end{bmatrix} \]
has a unique right minimal index, equal to $0$, and a unique nonzero partial multiplicity at $0$, equal to $1$. Any right minimal basis has the form 
\[ v=\alpha \begin{bmatrix}
1\\
-1
\end{bmatrix}, \qquad 0 \neq \alpha \in \C;\]
although the minimal basis is not unique, any has degree $0$ and thus encodes correctly the information on the minimal index. Similarly, it is not hard to check that any root polynomial has the form
\[r(\la)= \begin{bmatrix}
1+a(\la)+\la b(\la)\\
1-a(\la)+\la c(\la)
\end{bmatrix}, \qquad a(\la),b(\la),c(\la) \in \C[\la].\]  Any such root polynomial also forms a maximal set, as can be checked by the definition; in spite of the arbitrariness of the polynomials $a(\la),b(\la),c(\la)$, the order is always $1$ since \[ P(\la)r(\la)= (2\la + \la^2 b(\la) + \la^2 c(\la))  \begin{bmatrix}
1\\
1
\end{bmatrix},\]
and hence, in accordance to the theory, it corresponds to the partial multiplicity of the eigenvalue $0$. On the other hand, a generic zero direction may have an order that does not correspond to any minimal index or partial multiplicity. Indeed, let $k \in \mathbb{N}$ be any nonnegative integer, then
\[ z(\la) = \begin{bmatrix}
\la^k+1\\
\la^k-1
\end{bmatrix} \Rightarrow z(0)\neq 0, P(\la) z(\la) = 2\la^{k+1} \begin{bmatrix}
1\\
1
\end{bmatrix}\]
and thus $z(\la)$ is a zero direction of order $k+1$. In other words, in the case of singular polynomials the zero directions do not necessarily provide the correct information on the partial multiplicities.
\end{example}

The discussion above emphasizes that maximal sets of root polynomials are important bases that enclose the information on partial multiplicities, in a similar manner to how minimal bases enclose the information on minimal indices. (Although, unlike for minimal bases, the information on the partial multiplicity is not given by the degree but by the order.) A relevant question is therefore how to compute a maximal set of root polynomials, given a polynomial matrix $P(\la)$ and one point $\lambda_0$. If the calculation is performed numerically, it is also of interest to investigate the stability of any proposed algorithm.

A common approach to solve polynomial eigenvalue problems is to first linearize them: a pencil $L(\la)$ is called a linearization of $P(\la)$ if there exist unimodular (that is, invertible over $\C[\la]$) matrices $M(\la),N(\la)$ such that $M(\la) L(\la) N(\la) = P(\la) \oplus I$. If $L(\la)$ is a linearization of $P(\la)$, then all the finite zeros of $L(\la)$ and $P(\la)$ have the same nonzero partial multiplicities. In \cite[Section 8]{DopN21}, it was shown that for a very broad class of classical linearizations of polynomial matrices, including for example companion matrices, $\mathbb{L}_1$ and $\mathbb{L}_2$ vector spaces of linearizations, Fiedler pencils, and block Kronecker linearizations, it is very easy to recover a maximal set of root polynomials at $\lambda_0$ for $P(\la)$ from a maximal set of root polynomials at $\lambda_0$ for its linearization $L(\la)$: indeed, extracting a certain block suffices in all those cases. For more details, see in particular \cite[Theorem 8.5]{DopN21}, \cite[Theorem 8.9]{DopN21} and \cite[Theorem 8.10]{DopN21}.

For this reason, a robust algorithm for the computation of a maximal set of root polynomials at $\la_0$ for a pencil would immediately yield a robust algorithm for the computation of a maximal set of root polynomials at $\la_0$ for any polynomial matrix, namely:
\begin{enumerate}
\item Linearize $P(\la)$ via one of the linearizations for which recovery of maximal sets is described in \cite[Section 8]{DopN21}, say, $L(\la)$;
\item Compute a maximal set of root polynomials at $\la_0$ for $L(\la)$;
\item Extract a maximal set of root polynomials at $\la_0$ for $P(\la)$.
\end{enumerate}
This justifies a peculiar computational focus on the pencil case. The goal of this paper is to derive an algorithm for step 2 above;
our algorithm can in addition also compute minimal bases, which for many linearizations can also in turn be used to determine the minimal bases of the linearized polynomial matrix \cite{DDM09,DDM10,DLPVD,NP16}. In future work, we plan to investigate algorithms that work directly on the polynomial matrix $P(\la)$, and compare them with the approach described above.

\section{Zero directions of pencils} \label{Sec:Pencils}

As discussed above, the case of pencils is especially important because the existence of an algorithm to compute maximal sets of root polynomials and minimal bases for pencils immediately implies the existence of a general algorithm. We start in this section by considering the computation of zero directions: then, we will show how to extract root polynomials and minimal bases from them.

Finding the zero directions of an $m\times n$ pencil $L(\lambda) = L_0+\la L_1$ is a simpler problem than its analogue for a higher degree polynomial matrix,
since one can use the generalized Schur form of an arbitrary pencil which is related to the Kronecker canonical form.
Since the expansion of the pencil around $\la_0$ is again a pencil
$L(\lambda) = \hat{L}_0 + (\lambda-\lambda_0)\hat{L}_1,$
we can assume without loss of generality that the eigenvalue we are interested in is $\la_0=0$.

As a first step, we recall results from the literature concerning certain block Toeplitz matrices \cite{EliK98,KarK86,vdv}.
The first two allow us to retrieve the right minimal indices  $\epsilon_i, i=1,\ldots,n-r$ 
and the left minimal indices  $\eta_i, i=1,\ldots,m-r$ of the $m\times n$ pencil 
$L_0+\la L_1$ of normal rank $r$. 
\begin{theorem}[\cite{KarK86,Kar94}] \label{th:T1}
Let us denote the right nullity of the bidiagonal Toeplitz matrix
\begin{equation*}  \hat T_{k} := \left[ \begin{array}{cccc}
L_1 &  & \\
L_0 & \ddots &  \\
& \ddots & L_1 \\
& & L_0
\end{array}\right] \in \C^{m(k+1)\times nk},
\end{equation*}
by $\mu_k$, and set $\mu_0=\mu_{-1}=0$, and let $m_j$ be the number of 
right minimal indices  $\epsilon_i$ of the pencil $L_0+\la L_1$ equal to $j \geq 0$. Then
\[m_j=\mu_{j-1}-2\mu_j+\mu_{j+1}. \]
\end{theorem}

\begin{theorem}[\cite{KarK86,Kar94}] \label{th:T2}
Let us denote the left nullity of the bidiagonal Toeplitz matrix
\begin{equation*} \label{bidiag1} \tilde T_{k} := \left[ \begin{array}{ccccc}
L_0 & L_1 &  & \\
 & \ddots & \ddots & \\
 & & L_0 & L_1
\end{array}\right] \in \C^{mk\times n(k+1)},
\end{equation*}
by $\nu_k$, and set $\nu_0=\nu_{-1}=0$, and let $n_j$ be the number of 
left minimal indices  $\eta_i$ of the pencil $L_0+\la L_1$ equal to $j \geq 0$. Then
\[n_j=\nu_{j-1}-2\nu_j+\nu_{j+1}. \]
\end{theorem}

The third block Toeplitz result specializes the result mentioned earlier and relates 
to the elementary divisors at the eigenvalue $0$.
\begin{theorem}[\cite{vdv}] \label{th:T3}
Let us denote the rank of the bidiagonal Toeplitz matrix
\begin{equation*}  T_{k} := \left[ \begin{array}{ccccc}
L_0 & L_1 &  & \\
& L_0 & \ddots &  \\
& &  \ddots & L_1 \\
& & & L_0
\end{array}\right] \in \C^{m(k+1)\times n(k+1)},
\end{equation*}
by $r_k$, and set $r_{-1}=r_{-2}=0$, then the number of elementary divisors $\la^i$ of degree $i \geq 1$ is given by
$$ e_i=r_{i-2}-2r_{i-1}+r_{i}.
$$ 
\end{theorem}

\begin{remark}
Several papers \cite{Antoniou,EliK98,KarK86,VDooren93,vdv,Zuniga} have in the past made the link between the nullspaces of the block Toeplitz matrices mentioned in Theorems \ref{th:T1}, \ref{th:T2} and \ref{th:T3}. These earlier algorithms, however, do not directly address the problem of computing root polynomials, but they consider subproblems also mentioned in this paper. Comments on their complexity are given in Section \ref{Sec:Numerics}.
\end{remark}

One could in principle use the above Toeplitz matrices to compute the indices $m_i$, $n_i$ and $e_i$ and construct from
these minimal bases and root polynomials, but this would be very inefficient. The output of the
staircase algorithm \cite{vd79} applied to $L_0$ and $L_1$ in fact can be used to retrieve
all the information to find these polynomial vectors, as well as their degrees or orders. However, while the minimal indices and the partial multiplicities can be read out directly from the staircase form, the computation of a maximal set of root polynomials and a minimal basis requires some extra work, which is the subject of this paper; this is not dissimilar to what happens after having computed the Schur form of a matrix, from which the eigenvalues can be read directly while computing eigenvectors requires some further computational effort. 
There exists a staircase form for reconstructing both the left and right root polynomials and minimal bases, but since they are just the conjugate transpose of each other, we restrict ourselves here to the right case.

It is shown in
\cite{vd79} that there always exist unitary transformations
$U$ and $V$ (these will be real and orthogonal when the system and $\lambda_0$
are real) such that :
$$  U^\T (L_0+ \la L_1)V = A+ \la E:=
 \left[\begin{array}{c|c}
\tilde A  & \times \\  \hline 0 & A_{r} \end{array} \right] 
  + \la \left[\begin{array}{c|c}
\tilde E  & \\  \hline
0 & E_{r} \end{array} \right] $$
\begin{equation} \label{stair1} :=
\left[\begin{array}{cccc|c}
0 & A_{1,2} & \ldots & A_{1,k} &  \\
  &  \ddots &  \ddots & \vdots & \times \\
  &  & \ddots & A_{k-1,k} & \\
  &  &      & 0 & \\  \hline
  &  &      &   & A_{r} \end{array} \right] 
  + \la \left[\begin{array}{cccc|c}
E_{1,1} & E_{1,2} & \ldots  & E_{1,k}  &  \\
  &  \ddots &  \ddots & \vdots & \times \\
  &  & \ddots & E_{k-1,k} & \\
  &  &      & E_{k,k} & \\  \hline
  &  &      &   & E_{r} \end{array} \right]
\end{equation}

\noindent
where\\
(i) the matrices $E_{i,i}$ are of dimension $s_i \times t_i$ and of full
row rank $s_i$,\\
(ii) the matrices $A_{i,i+1}$ are of dimension $s_i \times t_{i+1}$ and
of full column rank $t_{i+1}$,\\
(iii) $A_r$ is of full column rank.\\
Hence, it follows \cite{vd79} that
\begin{equation*} t_1 \geq s_1 \geq t_2 \geq s_2 \geq \ldots \geq t_{k}
\geq s_{k} \geq 0 (:= t_{k+1}) 
\end{equation*}
and that the leading diagonal pencil has as structural elements associated with its local Smith form at $0$

\medskip

$m_i:=t_i-s_i$ right minimal indices equal to $i-1$, and 

$e_i:= s_i-t_{i+1}$ elementary divisors of degree $i$ at $0$.

\medskip

The transformations $U$ and $V$ are chosen to be unitary (or orthogonal in the real case) for reasons of numerical stability. 
The construction of the transformations is via the staircase algorithm, which recursively constructs growing othonormal bases $U_k$ and $V_k$ for the so-called Wong chains (or sequences) \cite{BerT12,Wong} defined as follows :
$$ \U_0=\{ 0 \}\in \C^m, \quad \V_i=L_0^{\leftarrow}\U_{i-1}\in \C^n, \quad \U_i= L_1\V_i\in \C^m, \quad \mathrm{for} \: i=1:k.
$$
Here we have borrowed from \cite{NP} the following notation to indicate how a matrix acts on a vector space:
$$  L_1\V:= \{u \;|\; u=L_1v, v\in \V \} \quad \mathrm{and} \quad L_0^{\leftarrow}\U:= \{ v \;|\; L_0 v=u, u\in \U \} .
$$
In other words, $L_1\V$ is the image of $\V$ under the transformation represented (in the canonical basis) by $L_1$ and $L_0^{\leftarrow}\U$ is the pre-image of $\U$ under the
transformation represented (in the canonical basis) by $L_0$.  Moreover, it is known \cite{BerT12,NP,vd79,Wong} that these spaces are nested and have the following dimensions
and {\em orthonormal} bases~:
\begin{equation} \label{Wong}  \U_k=\im U_k:= \im U\left[\begin{array}{c}
I_{\sigma_k} \\ 0 \end{array} \right], \; \sigma_k:=\sum_{i=1}^k s_i,
\quad 
  \V_k=\im V_k:= \im V\left[\begin{array}{c}
I_{\tau_k} \\ 0 \end{array} \right], \; \tau_k:=\sum_{i=1}^k t_i.
\end{equation}

\medskip

Moreover, since these
are constant transformations, they only transform the coordinate system in
which we have to construct the coefficients of the zero directions $x(\la)$ and
$y(\la)$. The above form \eqref{stair1} is appropriate for constructing
$x(\la)$, but there exists a dual form where the role of columns
and rows is interchanged, and which can thus be used to construct
$y(\la)$. Below we focus on finding solutions
$x_V(\la)$ in the coordinate system of \eqref{stair1}, but this is no loss of generality as the
corresponding zero directions of $L_0+\la L_1$ are easily seen to
be $x(\la)=V x_V(\la)$ (see also \cite[Proposition 3.2]{DopN21} for a more general result on root polynomials, and note that its proof can be adapted to show an analogous results on any zero direction).

\medskip

As a next step, in Section \ref{Sec:Extraction}, we will further simplify the pencil $A + \la E$, eventually reaching a form that allows us to (a) reduce the problem of computing a maximal set of root polynomials to the case of a regular pencil with only eigenvalues at $0$ and (b) reduce the problem of computing a minimal basis to the case of a pencil having only right minimal indices.

\section{Extracting the null space and root polynomials} \label{Sec:Extraction}

Although the Wong chains described in the previous section are uniquely defined as nested subspaces, the corresponding bases are not unique. It was shown in \cite{BeeV88} that an appropriate updating of the transformations $U$ and $V$ in \eqref{stair1}
makes sure that the ``stairs" $A_{i,i+1}$ and   $E_{i,i}$ have the following quasi-triangular form
\begin{equation} \label{stair2} A_{i,i+1}\Rightarrow \left[\begin{array}{c} \hat A_{i,i+1} \\ 0 \end{array}\right], \quad E_{i,i}\Rightarrow \left[\begin{array}{cc} 0 &  \hat E_{i,i} \end{array}\right],
\end{equation}
where $\hat A_{i,i+1}\in \C^{t_{i+1}\times t_{i+1}}$ and $\hat E_{i,i}\in \C^{s_i\times s_i}$ are both upper triangular and 
invertible. Such a form can be obtained by updating the transformations $U$ and $V$ to $\hat U=U U_d$ and $\hat V=V V_d$, using block-diagonal unitary matrices 
\begin{equation*} 
U_d:=\diag(I_{s_1},\hat U_{2,2}, \ldots, \hat U_{k,k}), \; \hat U_{i,i}\in \C^{s_i\times s_i},\quad V_d:=\diag(\hat V_{1,1}, \hat{V}_{2,2}, \ldots, \hat V_{k,k}), \; \hat V_{i,i}\in \C^{t_i\times t_i}.
 \end{equation*}
Note that these transformations have to be constructed backwards, starting with $\hat V_{k,k}$~:
$$  \left[\begin{array}{cc} 0 &  \hat E_{k,k} \end{array}\right] := E_{k,k}\hat V_{k,k}, $$
then for $i=k-1:-1:1$
$$ 
 \left[\begin{array}{c}  \hat A_{i,i+1} \\ 0 \end{array}\right] := \hat U_{i,i}(A_{i,i+1}\hat V_{i+1,i+1}), \quad
\left[\begin{array}{cc} 0 &  \hat E_{i,i} \end{array}\right] := (\hat U_{i,i}E_{i,i})\hat V_{i,i}.
$$
It is worth noting that the block columns of the updated transformations $\hat U$ and $\hat V$ are stil orthonormal bases for the nested Wong spaces defined earlier. We have only updated the choice of basis vectors. One can then separate the right null space blocks from the structure at the eigenvalue 0 by using the following result
\begin{lemma} \label{bidiag}
Let  $\tilde A+\lambda \tilde E$ be the leading principal subpencil of \eqref{stair1} and assume it is in staircase form with ``stairs" in the special form \eqref{stair2},
\begin{equation}\label{notbidiagonal}  \tilde A+ \la \tilde E:=
\left[\begin{array}{cccc}
0 & A_{1,2} & \ldots & A_{1,k} \\
  &  \ddots &  \ddots & \vdots  \\
  &  & \ddots & A_{k-1,k} \\
  &  &   & 0 
\end{array} \right] 
  + \la \left[\begin{array}{cccc}
E_{1,1} & E_{1,2} & \ldots  & E_{1,k}   \\
  &  \ddots &  \ddots & \vdots \\
  &  & \ddots & E_{k-1,k} \\
  &  &  & E_{k,k} \end{array} \right]
\end{equation}
indicating that this subpencil has only zero eigenvalues and right minimal indices. Then there exist unit upper triangular transformations $S$ and $T$ that eliminate all blocks except the 
 rank carrying stairs $E_{i,i}$ and $A_{i,i+1}$~:
\begin{equation} \label{bidiagonal} S^{-1}(\tilde A+\lambda \tilde E)T= \left[\begin{array}{cccc}
\la E_{1,1}  & A_{1,2} & 0 & 0 \\
  &  \ddots &  \ddots & 0  \\
  &  & \ddots & A_{k-1,k} \\ [1mm]
  &  &  & \la E_{k,k} 
\end{array} \right],
\end{equation}
without altering  $E_{i,i}$ and $A_{i,i+1}$, and hence puts the pencil in block-bidiagonal form.
\end{lemma}
\begin{proof}
The transformations $S$ and $T$ can be constructed recursively as a product of unit upper triangular transformations.
Indeed we can work upwards from block row $k-1$ till block row 1, each time eliminating first the elements in $\tilde A$ by row transformations, 
then the elements in $\tilde E$ by column transformations. The precise order of elimination is 
\begin{equation*}  
E_{k-1,k} \rightarrow  A_{k-2,k}  \rightarrow  [E_{k-2,k-1}, E_{k-2,k}]  \rightarrow  [A_{k-3,k-1}, A_{k-3,k}] \rightarrow \ldots  
\end{equation*}
The row transformations use full column rank matrices $A_{i-1,i}:={ \left[\begin{array}{cc}\hat A_{i-1,i} \\ 0 \end{array}\right]}$ above as pivot. 
The column transformations use the full row rank matrices $E_{i,i}:={ \left[\begin{array}{cc}0 & \hat E_{i,i} \end{array}\right]}$  to the left as pivot.
\end{proof}

In the proof above, the order in which the zero blocks of $\tilde A$ and $\tilde E$ are created is crucial in order to avoid destroying previously created zero blocks. For this reason it is necessary that this recurrence runs backwards.
A permuted version of this lemma was already shown in \cite{vd79}, without insisting that this only requires unit upper triangular transformation matrices. By choosing unit upper triangular matrices for this elimination, we can interpret it as the back-substitution for solving a linear system of equations in the unknowns $S_{i,j}$ and $T_{i,j}$~:
\begin{equation*}  \left(
\left[\begin{array}{cccc}
0 & A_{1,2} & \ldots & A_{1,k} \\
  &  \ddots &  \ddots & \vdots  \\
  &  & \ddots & A_{k-1,k} \\
  &  &   & 0 
\end{array} \right] 
  + \la \left[\begin{array}{cccc}
E_{1,1} & E_{1,2} & \ldots  & E_{1,k}   \\
  &  \ddots &  \ddots & \vdots \\
  &  & \ddots & E_{k-1,k} \\
  &  &  & E_{k,k} \end{array} \right] \right)  \left[\begin{array}{cccc}
I_{t_1} & T_{1,2} & \ldots  & T_{1,k}   \\
  &  \ddots &  \ddots & \vdots \\
  &  & \ddots & T_{k-1,k} \\
  &  &  & I_{t_k} \end{array} \right]  \end{equation*}
  \begin{equation}\label{ecc}
  = \left[\begin{array}{ccccc}
I_{s_1} & S_{1,2} & \ldots  & S_{1,k-1} & 0   \\
  &  \ddots &  \ddots & \vdots & \vdots\\
  &  & \ddots & S_{k-2,k-1} & 0 \\
  &  &  & I_{s_{k-1}}  & 0 \\
  &  &  &  & I_{s_k} \end{array} \right]  \left[\begin{array}{cccc}
\la E_{1,1}  & A_{1,2} & 0 & 0 \\
  &  \ddots &  \ddots & 0  \\
  &  & \ddots & A_{k-1,k} \\ [1mm]
  &  &  & \la E_{k,k} 
\end{array} \right].
\end{equation}
Let us write the matrices $\tilde A$, $\tilde E$, $S$ and $T$ as follows
$$ \tilde A = A_{d}+A_u, \quad \tilde E = E_{d}+E_u, \quad S = I_n + S_u, \quad T= I_m + T_u,  
$$
where $A_d$ and $E_d$ are the submatrices of $\tilde A$ and $\tilde E$ that are kept in the bidiagonal pencil.
Then the equations \eqref{ecc} can be rewritten as
$$ S_uA_d-(A_d+A_u)T_u= A_u, \quad   S_uE_d-(E_d+E_u)T_u = E_u
$$
where  the submatrices 
$$ A_{i,j}, 1\le i \le k-2,  \; i+1<j\le k, \qquad E_{i,j}, \; 1\le i \le k-1,  \; i<j\le k$$
are to be eliminated form the right hand sides, and the submatrices 
$$ 
S_{i,j}:={\small \left[\begin{array}{cc} \hat S_{i,j} & 0 \end{array}\right]}, \; 1\le i \le k-2,  \; i<j\le k-1, \quad T_{i,j}:={\small \left[\begin{array}{cc} 0 \\ \hat T_{i,j} \end{array}\right]},  \; 1\le i \le k-1,  \; i<j\le k 
$$ 
are the unknowns and have as many nonzeros as the blocks they are supposed to eliminate. Their block structure follows from the block structure of the pivot blocks, as indicated in the proof of Lemma \ref{bidiag}.
This system of equations is therefore invertible and we can apply iterative refinement \cite{Higham} to improve the accuracy of the corresponding computations. Since there is no factorization to be performed for the iterative refinement, 
its computational cost is very reasonable. Also the choice of transformation is such that $|\det \hat V T|=|\det \hat U S|=1$.

Notice that the Wong sequences are now also spanned by the growing subblocks of the transformation matrices $\hat U S$ and $\hat V T$ but these bases are no longer orthonormal.  
\begin{corollary}
The matrices $\hat S:=\hat US=U U_d S$ and $\hat T:=\hat V T=VV_d T$ have non orthonormal block columns that still span the nested Wong spaces defined in \eqref{Wong}
$$  \U_k = \im \hat S \left[\begin{array}{c}
I_{\sigma_k} \\ 0 \end{array} \right], \; \sigma_k:=\sum_{i=1}^k s_i,
\quad 
  \V_k=\im \hat T \left[\begin{array}{c}
I_{\tau_k} \\ 0 \end{array} \right], \; \tau_k:=\sum_{i=1}^k t_i.
$$
\end{corollary}
We point out that the staircase form index sets $\{s_i\}$ and $\{t_i\}$ allow us to separate the pencil in its
two distinct structures, namely the right minimal indices and the Jordan structure at the eigenvalue $0$.

\begin{example} \label{matlabex}
We give an example of a special case of \eqref{notbidiagonal} illustrating this separation for the dimensions 
$$ k=3, \quad t_1=5, \; s_1=4,\;  t_2=3,\;  s_2=2,\; t_3=1,\; s_3=0,\; t_4=0.
$$
\begin{equation*} \left[\begin{array}{ccc:cc|cc:c|c}
 \zr & \zr & \zr & \zg & \zg & \otr & \tmr & \tmg & \tmr \\
 \zr & \zr & \zr & \zg & \zg & \zr & \otr & \tmg & \tmr \\ \hdashline
\zg & \zg & \zg & \zb & \zb & \zg & \zg & \otb & \tmg \\ 
\zg & \zg & \zg & \zb & \zb & \zg & \zg & \zb & \tmg \\ \hline
\zr & \zr & \zr & \zg & \zg & \zr & \zr & \zg & \otr \\  \hdashline
\zg & \zg & \zg & \zb & \zb & \zg & \zg & \zb & \zg
 \end{array}\right] + \la
 \left[\begin{array}{ccc:cc|cc:c|c}  \zr & \otr & \tmr & \tmg & \tmg & \tmr & \tmr & \tmg & \tmr \\
 \zr & \zr & \otr & \tmg & \tmg & \tmr & \tmr & \tmg & \tmr \\ \hdashline
\zg & \zg & \zg & \otb & \tm & \tmg & \tmg & \tm & \tmg \\ 
\zg & \zg & \zg & \zb & \otb & \tmg & \tmg & \tm & \tmg \\ \hline
\zr & \zr & \zr & \zg & \zg & \zr & \otr & \tmg & \tmr \\  \hdashline
\zg & \zg & \zg & \zb & \zb & \zg & \zg & \otb & \tmg  \end{array}\right] 
\end{equation*}
The red blocks correspond to the right minimal indices and the blue blocks correspond to the Jordan structure at 0; the symbol $\times$ denotes an arbitrary complex number that is allowed to be nonzero while the symbol $\otimes$ denotes a complex number that is guaranteed to be nonzero.
The sequence of indices for these two structures are 
$$ \color{red} k^{(r)}=3, \quad t_1^{(r)}=3, \; s_1^{(r)}=2, \; t_2^{(r)}=2, \; s_2^{(r)}=1,\; t_3^{(r)}=1,\; s_3^{(r)}
=0,\; t_4^{(r)}=0,
$$
and 
$$ \color{blue} k^{(b)}=3, \quad t_1^{(b)}=2, \; s_1^{(b)}=2, \; t_2^{(b)}=1, \; s_2^{(b)}=1,\; t_3^{(b)}=0,\; s_3^{(b)}=0,\; t_4^{(b)}=0,
$$
and are obtained by starting from ${\color{blue}t_4^{(b)}}=t_4$ and then equating (for decreasing $i$) 
$$  {\color{red}s_{i-1}^{(r)}}={\color{red}t_i^{(r)}}=t_i-{\color{blue}t_i^{(b)}}, \quad {\color{blue}t_{i-1}^{(b)}}={\color{blue}s_{i-1}^{(b)}}=s_{i-1}-{\color{red}s_{i-1}^{(r)}}, \quad \mathrm{for} \; i=k:-1:2, \quad {\color{red} t_1^{(r)}}=t_1-{\color{blue}t_1^{(b)}}.
$$
Unfortunately, these blocks are not completely decoupled in this coordinate system. This can be cured by the following reordering of rows and columns 
$$ [{\color{red}1,2},{\color{blue}3,4},{\color{red}5},{\color{blue}6}] 
\longrightarrow [{\color{red}1,2,5},{\color{blue}3,4,6}]   \quad 
\mathrm{and}  \quad  [{\color{red}1,2,3},{\color{blue}4,5},{\color{red}6,7},{\color{blue}8},{\color{red}9}] \longrightarrow [{\color{red}1,2,3,6,7,9},{\color{blue}4,5,8}] 
$$
(which is equivalent to multiplying by permutation matrices on the left and on the right) to separate the (red) Kronecker block structure from the (blue) Jordan structure.
\begin{equation*} \left[\begin{array}{cccccc|ccc}
\zr & \zr & \zr & \otr & \tmr & \tmr & \zg & \zg & \tmg \\
\zr & \zr & \zr & \zr & \otr & \tmr & \zg & \zg & \tmg \\
\zr & \zr & \zr & \zr & \zr & \otr & \zg & \zg & \zg \\ \hline
\zg & \zg & \zg & \zg & \zg & \tmg & \zb & \zb & \otb \\ 
\zg & \zg & \zg & \zg & \zg & \tmg & \zb & \zb  & \zb \\ 
\zg & \zg & \zg & \zg & \zg & \zg & \zb & \zb & \zb 
\end{array}\right] + \la
\left[\begin{array}{cccccc|ccc} 
\zr & \otr & \tmr & \tmr & \tmr & \tmr & \tmg & \tmg & \tmg  \\
\zr & \zr & \otr  & \tmr & \tmr & \tmr & \tmg & \tmg & \tmg \\
\zr & \zr & \zr  & \zr & \otr & \tmr & \zg & \zg & \tmg \\ \hline
\zg & \zg & \zg  & \tmg & \tmg & \tmg & \otb & \tm & \tm \\ 
\zg & \zg & \zg  & \tmg & \tmg & \tmg & \zb & \otb & \tm \\ 
\zg & \zg & \zg & \zg & \zg  & \tmg & \zb & \zb & \otb  \end{array}\right].
\end{equation*}
Moreover, applying the same reordering on the bidiagonal pencil described in Lemma \ref{bidiag} would yield the required block triangular form. The bidiagonal form for our example is
\begin{equation*} \left[\begin{array}{ccc:cc|cc:c|c}
 \zr & \zr & \zr & \zg & \zg & \otr & \tmr & \tmg & \zr \\
 \zr & \zr & \zr & \zg & \zg & \zr & \otr & \tmg & \zr \\  \hdashline
\zg & \zg & \zg & \zb & \zb & \zg & \zg & \otb & \zg \\ 
\zg & \zg & \zg & \zb & \zb & \zg & \zg & \zb & \zg \\ \hline
\zr & \zr & \zr & \zg & \zg & \zr & \zr & \zg & \otr \\  \hdashline
\zg & \zg & \zg & \zb & \zb & \zg & \zg & \zb & \zg
 \end{array}\right] + \la
 \left[\begin{array}{ccc:cc|cc:c|c}  \zr & \otr & \tmr & \tmg & \tmg & \zr & \zr & \zg & \zr \\
 \zr & \zr & \otr & \tmg & \tmg & \zr & \zr & \zg & \zr \\  \hdashline 
\zg & \zg & \zg & \otb & \tm & \zg & \zg & \zb & \zg \\ 
\zg & \zg & \zg & \zb & \otb & \zg & \zg & \zb & \zg \\ \hline
\zr & \zr & \zr & \zg & \zg & \zr & \otr & \tmg & \zr \\  \hdashline
\zg & \zg & \zg & \zb & \zb & \zg & \zg & \otb & \zg  \end{array}\right] 
\end{equation*}
and its permuted form is block upper triangular
\begin{equation*} \left[\begin{array}{cccccc|ccc}
\zr & \zr & \zr & \otr & \tmr & \zr & \zg & \zg & \tmg \\
\zr & \zr & \zr & \zr & \otr & \zr & \zg & \zg & \tmg \\
\zr & \zr & \zr & \zr & \zr & \otr & \zg & \zg & \zg \\ \hline
\zg & \zg & \zg & \zg & \zg & \zg & \zb & \zb & \otb \\ 
\zg & \zg & \zg & \zg & \zg & \zg & \zb & \zb  & \zb \\ 
\zg & \zg & \zg & \zg & \zg & \zg & \zb & \zb & \zb 
\end{array}\right] + \la
\left[\begin{array}{cccccc|ccc} 
\zr & \otr & \tmr & \zr & \zr & \zr & \tmg & \tmg & \zg  \\
\zr & \zr & \otr  & \zr & \zr & \zr & \tmg & \tmg & \zg \\
\zr & \zr & \zr  & \zr & \otr & \zr & \zg & \zg & \tmg \\ \hline
\zg & \zg & \zg  & \zg & \zg & \zg & \otb & \tm & \zb \\ 
\zg & \zg & \zg  & \zg & \zg & \zg & \zb & \otb & \zb \\ 
\zg & \zg & \zg & \zg & \zg  & \zg & \zb & \zb & \otb  \end{array}\right] .
\end{equation*}

\end{example}

\medskip

Equivalently, a complete block diagonal decoupling is also obtained if we update the unit upper triangular matrices $S$ and $T$ to block diagonalize the upper-triangular matrices $\hat A_{i-1,i}$ and $\hat E_{i,i}$ in the bidiagonal form of Lemma \ref{bidiag}. 

\section{Recurrences for the null space and root polynomials} \label{Sec:Recurrences}

The discussion in Section \ref{Sec:Extraction} shows that we can obtain a triangular block decomposition of the form
\begin{equation}\label{eq:goodform}
\hat S^{-1}(L_0+\la L_1)\hat T =\left[\begin{array}{ccc}
  A_\epsilon+\la E_\epsilon & 0   & \times \\
0  & A_0+\la E_0  & \times \\  0  &  0 & A_{r}+\la E_r \end{array} \right],
\end{equation} 
where $A_\epsilon+\la E_\epsilon$ only has a right null space structure, $A_0+\la E_0 $ only has  a Jordan structure at $\la=0$, and $A_r + \la E_r$ has full column rank and contains the rest of the pencil structure. Moreover, $A_\epsilon+\la E_\epsilon$ and $A_0+\la E_0$ are both in a bidiagonal staircase form given in \eqref{stair2} and \eqref{bidiagonal}, and the equivalence transformation pair $(\hat S,\hat T)$ was obtained as the product $(\hat U \cdot S \cdot \Pi_1,\hat V\cdot T \cdot \Pi_2)$ of a unitary equivalence transformation pair $(\hat U,\hat V)$, an upper triangular equivalence transformation pair $(S,T)$, and a permutation transformation pair $(\Pi_1, \Pi_2)$. This decomposition could be important for future research on the analysis of
the numerical stability of each step of the method since the unitary similarity was analyzed in \cite{vd79}, and the upper triangular 
similarity pair can be viewed as a back-substitution step for solving a linear system, whose accuracy can be improved using iterative refinement \cite{Higham}. However, a detailed analysis of the composed algorithm is beyond the scope of the present paper.
 
\subsection{Calculating root polynomials: reduction to the regular case}

We now argue that the structured pencils \eqref{eq:goodform} allows us to just compute the root polynomials of the pencil $A_0 + \la E_0$, which is guaranteed to be regular and only have Jordan blocks with eigenvalue $0$. We start by defining what it means for a polynomial matrix to be right invertible.

\begin{definition}
We say that $A(\la) \in \C[x]^{m \times n}$ is right invertible over $\C[x]$ if there exists a polynomial $A(\la)^R \in \C[x]^{n \times m}$ such that $A(\la) A(\la)^R=I_m$. If $m=0$, by convention the $0 \times n$ empty matrix polynomial is right invertible and its right inverse is its $n \times 0$ transpose.
\end{definition}

Linear equations whose cofficients are right invertible polynomial matrices always have polynomial solutions in the sense that if $A(\la)$ is right invertible with right inverse $A(\la)^R$ then the $\C[\la]$-linear map associated with $A(\la)$ is surjective and thus the equation $A(\la) y(\la) = b(\la)$, $b(\la) \in \C[\la]^m$, has at least one polynomial solution. Indeed, a solution can be constructed as $y(\la)=A(\la)^R b(\la) \in \C[\la]^n$. Lemma \ref{lem:rightinvertible} characterizes completely the set of right invertible polynomial matrices in terms of their eigenvalues and minimal indices. We note that almost equivalent results had appeared in \cite{BeeVsimax}; however, Lemma \ref{lem:rightinvertible} below provides a slightly better bound on the degree of a right inverse, and for this reason we will include a proof. 


\begin{lemma}\label{lem:rightinvertible}
Let $L(\la) \in \C[\la]^{m \times n}$ be a matrix pencil. Then, $L(\la)$ is right invertible over $\C[\la]$ if and only if it has no finite eigenvalues and no left minimal indices. Moreover, in that case, there exists a right inverse having degree $\leq m-1$.
\end{lemma}
\begin{proof}
First, we point out that \begin{enumerate}
 \item for any pair of square invertible matrices $S,T$ it holds that $L(\la)$ is right invertible with right inverse $R(\la)$ if and only if $SL(\la)T$ is right invertible with right inverse $T^{-1} R(\la) S^{-1}$;
 \item $L_1(\la),L_2(\la)$ are both right invertible with right inverses $R_1(\la),R_2(\la)$ resp. if and only if $L_1(\la) \oplus L_2(\la)$ is right invertible with right inverse $R_1(\la) \oplus R_2(\la)$;
 \item The degree of a direct sum of polynomial matrices is the maximum of their degrees.
\end{enumerate}
Hence, we can assume without loss of generality that $L(\la)$ is a single block within a Kronecker canonical form.

Assume first that $L(\la)$ is either a left Kronecker block or a Jordan block associated with a finite eigenvalue. Then there exists a nonzero vector $u \in \C^m$ and a scalar $\mu \in \C$ such that $u^* L(\mu)=0$. Indeed, for a Kronecker block $\mu$ can be any element of $\C$ and $u$ can be constructed by evaluating a row of its left minimal basis at $\la=\mu$; while if $L(\la)$ is a Jordan block with finite eigenvalue $\la_0$ then we take $\mu=\la_0$  and $u^*$ equal to a left eigenvector. Assuming for a contradiction that $L(\la)$ is right invertible we then have for any right inverse $L(\la)^R$
\[ 0 \neq u^* = u^* L(\la) L(\la)^R \Rightarrow 0 \neq u^* L(\mu) L(\mu)^R = 0.\]

We now give a constructive proof of the converse implication, i.e., we explicitly compute a right inverse that also satisfies the degree bound. This time, we can assume that $L(\la)$ is either an $m \times (m+1)$ Kronecker block or an $m \times m$ Jordan block with infinite eigenvalues. We treat each case separately.
\begin{itemize}
\item \emph{Right singular Kronecker block.}  If $m=0$, the statement is true by definition. Otherwise,
\[ L(\la) = \begin{bmatrix}
I_m & 0_{m \times 1}
\end{bmatrix} + \la \begin{bmatrix}
0_{m \times 1}&  I_m
\end{bmatrix} \in \C[\la]^{m \times (m+1)}\]
and it suffices to take the Toeplitz right inverse
\[ L(\la)^R = \begin{bmatrix}
1&-\la&\la^2&\dots&(-\la)^{m-1}\\
0&1 &-\la&\dots&(-\la)^{m-2}\\
\vdots&\ddots&\ddots&\ddots&\vdots\\
0&\dots&0&1&-\la\\
0&\dots&&0&1\\
0&&\dots&&0
\end{bmatrix} \in \C[\la]^{(m+1) \times m}.\]
\item \emph{Jordan block at infinity.}
In this case $L(\la)=I_m + \lambda J$ where $J$ is a nilpotent $m \times m$ Jordan block and we can take the Toeplitz right inverse
\[ L(\la)^R = \begin{bmatrix}
1&-\la&\la^2&\dots&(-\la)^{m-1}\\
0&1 &-\la&\dots&(-\la)^{m-2}\\
\vdots&\ddots&\ddots&\ddots&\vdots\\
0&\dots&0&1&-\la\\
0&\dots&&0&1
\end{bmatrix} \in \C[\la]^{m \times m}.\]
\end{itemize}

\end{proof}

From now on we will assume that the pencil $L(\la)$ is in the form
\begin{equation}\label{eq:staircase}
L(\la)=\begin{bmatrix}
L_\ell(\la) & A(\la) & B(\la)\\
0 & R(\la) & C(\la)\\
0 & 0 & L_r(\la)
\end{bmatrix}
\end{equation} 
where $L_\ell(\la)$ is right invertible over $\C[\la]$, $L_r(\la)$ has no eigenvalues at $\la_0 \in \C$ and has no right minimal indices, and $R(\la)$ is regular and has only eigenvalues at $\la_0$. Note also that what we had achieved in \eqref{eq:goodform} is a special case of \eqref{eq:staircase} with $A(\la)=0$.

\begin{proposition}\label{prop:minb}
Suppose that $L(\la)$ is as in \eqref{eq:staircase} and satisfies the assumptions stated immediately below it. Then, the columns of $M(\la)$ are a minimal basis for $L_\ell(\la)$ if and only if a minimal basis for $L(\la)$ is given by the columns of
\[ \begin{bmatrix}
M(\la)\\
0\\
0
\end{bmatrix}.\]
\end{proposition}
\begin{proof}
Since $R(\la)$ and $L_r(\la)$ have both trivial right nullspaces, it is clear that any minimal basis for $L(\la)$ can only have nonzero top blocks. After this observation, the proof becomes trivial.
\end{proof}

\begin{proposition}\label{prop:extractrootpoly}
Suppose that $L(\la)$ is as in \eqref{eq:staircase} and satisfies the assumptions stated immediately below it. Then:
\begin{enumerate}
    \item  If $v(\la) + (\la-\lambda_0)^\ell z(\la)$, where $\deg v(\la)\leq \ell-1$, is a root polynomial of order $\ell$  for $R(\la)$ at $\lambda_0$ then there is a polynomial vector $u(\la)$ such that 
\[ q(\la)= \begin{bmatrix}
u(\la)\\
v(\la)\\
0
\end{bmatrix}\]
is a root polynomial of order exactly $\ell$ at $\lambda_0$  for $L(\la)$;
\item If \[ g(\la) = \begin{bmatrix}
a(\la)\\
b(\la)\\
c(\la)
\end{bmatrix} + (\la-\la_0)^\ell r(\la) \quad (\deg \begin{bmatrix}
a(\la)\\
b(\la)\\
c(\la)
\end{bmatrix}  \leq \ell-1) \] is a root vector of order $\ell$ at $\la_0$ for $L(\la)$, then $c(\la)=0$ and $b(\la)$ is a root vector of order exactly $\ell$ at $\la_0$ for $R(\la)$.
\end{enumerate}

\end{proposition}
\begin{proof}
\begin{enumerate}
\item Suppose that $R(\la)v(\la)=(\la-\lambda_0)^\ell w(\la)$ with $v(\lambda_0)\neq0\neq w(\lambda_0)$.  Then, by Proposition \ref{prop:minb}, $q(\la_0)\not \in \ker_{\la_0} L(\la)$ \cite[Definition 2.7]{DopN21} (see also \cite[Sec. 8]{N11}). Moreover,
\[ L(\la) b(\la) = \begin{bmatrix}
L_\ell(\la) u(\la) + A(\la) v(\la)\\
R(\la) v(\la)\\
0
\end{bmatrix}.\]
Hence, $q(\la)$ is a root polynomial of order exactly $\ell$ for $L(\la)$ provided that $L_\ell(\la) u(\la) + A(\la) v(\la) \equiv 0 \mod (\la-\la_0)^\ell$; but it follows by Lemma \ref{lem:rightinvertible} that choosing $u(\la)=-L_\ell(\la)^RA(\la)v(\la)$, where $L_\ell(\la)^R$ denotes any polynomial right inverse of $L_\ell(\la)$, suffices.
\item Now suppose that $g(\la)$ is a root polynomial of order $\ell$ at $\la_0$ for $L(\la)$, implying that $\begin{bmatrix}
a(\la)^T & b(\la)^T & c(\la)^T
\end{bmatrix}^T$ also is a root polynomial for $L(\la)$, of order $\geq \ell$.
Then,
\[ \begin{cases} L_\ell(\la) a(\la) + A(\la) b(\la) + B(\la)c(\la) \equiv 0 \mod (\la-\lambda_0)^\ell\\
R(\la) b(\la) + C(\la) c(\la) \equiv 0 \mod (\la-\lambda_0)^\ell \\
L_r(\la) c(\la) \equiv 0 \mod (\la-\lambda_0)^\ell
\end{cases}
\]
Writing $L_r(\la)=L_0 + L_1 (\la-\lambda_0)$ and $c(\la) \equiv \sum_{k=0}^{\ell-1}c_k(\la-\lambda_0)^k \mod (\la-\lambda_0)^\ell$, the last equation is equivalent to
\[ \begin{bmatrix}
L_0 & L_1 & &\\
& L_0 & \ddots & \\
& &\ddots & L_1 \\
& & &L_0
\end{bmatrix} \begin{bmatrix}
c_{\ell-1}\\
\vdots\\
c_1\\
c_0
\end{bmatrix} = 0.  \]
However, since $L_r(\la)$ has no eigenvalues at $\lambda_0$ and no right minimal index,  then $L_0=L_r(\lambda_0)$ must be left invertible over $\C$, and hence so is the coefficient matrix above (proof of the latter claim: the rank of a block triangular matrix is bounded below by the sum of the ranks of the diagonal blocks). We deduce that the associated $\C$-linear map is injective and thus $c(\la) \equiv 0 \mod (\la-\lambda_0)^\ell \Rightarrow c(\la)=0$ (since $c(\la)$ has degree $\ell-1$ at most). Hence, $R(\la)b(\la)$ must be a multiple of $(\la-\la_0)^\ell$.  Now, suppose for a contradiction that $b(\la_0) = 0$. Then, by Proposition \ref{prop:minb}, $a(\la_0) \not \in \ker_{\la_0} L_{\ell}(\la) $. On the other hand, $L_\ell(\la_0)a(\la_0)=-A(\la_0)b(\la_0)=0$ implying that $L_\ell(\la) a(\la) \equiv 0 \mod (\la-\la_0)$. Hence, $a(\la)$ is a root polynomial for $L_\ell(\la)$ at $\la_0$: a contradiction, as $L_\ell(\la)$ is right invertible and therefore by Lemma \ref{lem:rightinvertible} does not have finite eigenvalues. We conclude that $b(\la_0)\neq 0$ and therefore $b(\la)$ is a root polynomial of order $\geq \ell$ at $\la_0$ for $R(\la)$. Suppose now for a contradiction that the order is strictly greater than $\ell$; then, $\ell < \deg R(\la) b(\la) \leq \deg R(\la) + \deg b(\la)  \leq 1 + \ell - 1 = \ell$, which is absurd.
\end{enumerate} 

\end{proof}

\begin{remark}
Observe that, although we have proved Proposition \ref{prop:extractrootpoly} in a slightly more general setting (for the benefit of any reader who may wish to use the result under slightly more relaxed assumptions), our algorithmic construction leads to \eqref{eq:goodform} which is a special case of \eqref{eq:staircase} with $A(\la)=0$. The proof of  Proposition \ref{prop:extractrootpoly} makes it clear that, in that case, one can construct a root polynomial for $L(\la)$ from one of $R(\la)$ taking $u(\la)=0$.
\end{remark}

\begin{theorem}\label{thm:extractrootpolys}
In the notation of \eqref{eq:staircase} and under the assumptions of this subsection, $\{v_i(\la)\}_{i=1}^s$ is a maximal set of root polynomials  for $R(\la)$ at $\lambda_0$ of order $\ell_1, \dots \ell_s$ if and only if a maximal set of root polynomials (of the same orders)  for $L(\la)$ has the form $\{v_i(\la)\}_{i=1}^s$ where
\[ q_i(\la)=\begin{bmatrix}
u_i(\la)\\
v_i(\la)\\
0
\end{bmatrix} + (\la-\lambda_0)^{\ell_i} r_i(\la), \qquad i=1,\dots,s.\]
\end{theorem}

\begin{proof}
 Let $\{v_i(\la)\}_{i=1}^s$ be a maximal set of root polynomials at $\la_0$ for $R(\la)$, of orders $\ell_i$. It follows from Proposition \ref{prop:extractrootpoly} that we can construct a set of $\{ q_i(\la) \}_{i=1}^s$ of the sought form, each of which is a root polynomial of order $\ell_i$ for $L(\la)$. That the set $\{ q_i(\la) \}_{i=1}^s$ is $\la_0$-independent follows by Proposition \ref{prop:minb} and because, if $M(\la)$ has $p$ columns and is a minimal basis for $\ker L_\ell(\la)$, then
\[ p+s \geq \rank  \begin{bmatrix}
M(\la_0) & u_1(\la_0) & \dots & u_s(\la_0)\\
0 & v_1(\la_0) & \dots & v_s(\la_0)\\
0 & 0 & \dots & 0
\end{bmatrix} \geq \rank M(\la_0)+ \rank \begin{bmatrix}
v_1(\la_0) & \dots & v_s(\la_0)
\end{bmatrix}=p+s.\]
Completeness and maximality then follow from \cite[Theorem 3.10]{nofvandue}, Theorem \ref{thm:frovan} and the fact that the nonzero partial multiplicities associated with the eigenvalue $\lambda_0$ are the same for $L(\la)$ and $R(\la)$.

Conversely suppose that $\{ q_i (\la) \}_{i=1}^s$ is a maximal set of root polynomials at $\la_0$ for $L(\la)$, of orders $\ell_i$; Proposition \ref{prop:extractrootpoly} guarantees the bottom block of $q_i(\la)$ must be $0 \mod (\la-\la_0)^{\ell_i}$. Moreover, for each $i$, again by Proposition \ref{prop:extractrootpoly} the middle block $v_i(\la)$ is a root polynomial of order $\ell_i$ at $\la_0$ for $R(\la)$. Assume for a contradiction that $\{ v_i(\la)\}_{i=1}^s$ are not a $\la_0$-independent set: then there are coefficients $d_i$ not all zero and such that $ \sum_{i=1}^s d_i v_i(\la_0)=0 \Rightarrow w:= \sum_{i=1}^s d_i u_i(\la_0) \not \in \ker_{\la_0} L_\ell(\la)$
where the last implication follows from the $\la_0$-independence of the $q_i(\la)$. On the other hand,
$L_\ell(\la_0) w = - A(\la_0) \sum_{i=1}^s d_i v_i(\la_0)=0$
and hence $w$ is a root polynomial at $\la_0$ for $L_\ell(\la)$, contradicting Lemma \ref{lem:rightinvertible}. At this point, completeness and maximality follow by the same argument as above.

\end{proof}

\begin{remark}
Again, if $L(\la)$ is in the form \eqref{eq:goodform}, then one may take $u_i(\la)=0$ in Theorem \ref{thm:extractrootpolys}when constructing a maximal set for $L(\la)$ from a maximal set for $R(\la)$.
\end{remark}

\subsection{Constructing a (right) minimal basis}

In this subsection we focus on the calculation of a right minimal basis for an $m\times n$ pencil $A_{\epsilon}+\la E_{\epsilon}$ with only right minimal indices
$\{\epsilon_1,\ldots,\epsilon_{n-r}\}$, and in bidiagonal
staircase form
\begin{equation}\label{eq:ventisei}
A_\epsilon + \la E_\epsilon =
\left[\begin{array}{cccc}
 \la  E_{1,1} & A_{1,2} &   &  \\
 &  \ddots &  \ddots &  \\
 &  & \la E_{k-1,k-1}  & A_{k-1,k} 
\end{array} \right] ,
\end{equation}
where 
$$ A_{i,i+1} \in \C^{s_i\times s_{i}}, \quad  E_{i,i} = \left[\begin{array}{cc}0 & \hat E_{i,i}\end{array}\right]\in \C^{s_i\times t_i}, \quad \hat E_{i,i}\in \C^{s_i\times s_i}, 
$$
and both $A_{i,i+1}$ and $\hat E_{i,i}$ are invertible upper triangular matrices. Note that this implies that $t_{i+1}=s_i$ for $i=1:k-1$ and $s_k=t_{k+1}=0$. Moreover, note that by Proposition \ref{prop:minb} the task of computing a minimal basis for a pencil having the structure \eqref{eq:ventisei} is sufficient to compute a minimal basis for a pencil with the structure of \eqref{eq:goodform}, and therefore also for a general pencil by the previous analyses.

Let $N(\la)$ be a matrix whose columns are such a minimal basis. As pointed out earlier, the coefficients $N_i$ in $N(\la):= \sum_{i=0}^{k-1} N_i \la^{i}$ satisfy the convolution 
equation
$$ \left[ \begin{array}{cccc}
E_\epsilon &  & \\
A_\epsilon & \ddots &  \\
& \ddots & E_\epsilon \\
& & A_\epsilon
\end{array}\right] \left[ \begin{array}{c}
N_{k-1}\\
\vdots \\
N_1 \\
N_0
\end{array}\right] 
=0.
$$
Yet, it is more efficient to exploit the bidiagonal staircase form since
the matrix pencil \eqref{eq:goodform} has a right null space basis $N(\la)$ of the following form
\begin{equation}\label{eq:minimalbasis}
N(\la):=\left[ \begin{array}{c} I_{t_1} \\ Z_1 \la \\ \vdots \\ Z_{k-1}...Z_{1} \la^{k-1} \end{array}\right] \in \C[\la]^{n\times(n-r)}
\end{equation} 
where $n=\sum_{i=1}^k t_i$ and $n-r=\sum_{i=1}^k (t_i-s_i)= t_1$ and $Z_i$ is the minimum norm solution of
\begin{equation} \label{Zi}
 A_{i,i+1}Z_i + E_{i,i} =0, \quad \Rightarrow \; Z_i=\left[\begin{array}{cc} 0 & - A_{i,i+1}^{-1}\hat E_{i,i}\end{array}\right] \in\C^{s_{i}\times t_i}.
\end{equation}

\begin{theorem} \label{th:rightinv}
Let 
$$
A_\epsilon + \la E_\epsilon =
\left[\begin{array}{cccc}
 \la  E_{1,1} & A_{1,2} &   &  \\
 &  \ddots &  \ddots &  \\
 &  & \la E_{k-1,k-1}  & A_{k-1,k} 
\end{array} \right] , 
$$
be a bidiagonal pencil where 
$$ A_{i,i+1} \in \C^{s_i\times s_{i}}, \quad  E_{i,i} = \left[\begin{array}{cc}0 & \hat E_{i,i}\end{array}\right]\in \C^{s_i\times t_i}, \quad \hat E_{i,i}\in \C^{s_i\times s_i}, 
$$
and both $A_{i,i+1}$ and $\hat E_{i,i}$ are invertible and upper triangular.
Then the columns of
$$N(\la):=\left[ \begin{array}{c} I_{t_1} \\ Z_1 \la \\ \vdots \\ Z_{k-1}...Z_{1} \la^{k-1} \end{array}\right] \in \C[\la]^{n\times(n-r)}\quad \mathrm{where} \quad Z_i=\left[\begin{array}{cc} 0 & - A_{i,i+1}^{-1}\hat E_{i,i}\end{array}\right] \in\C^{s_{i}\times t_i}
$$
form a minimal polynomial basis for its right null space.
\end{theorem}
\begin{proof}
It suffices \cite{For75} to show that (a) $N(\mu)$ has full rank for every $\mu \in \C$ and (b) $N(\la)$ is column reduced. (a) is trivial since the top $t_1 \times t_1$ submatrix of $N(\mu)$ is always the identity, for all $\mu$. For (b), observe that for all $h=1,\dots,k-1$ it holds \[\prod_{i=1}^h Z_i = \begin{bmatrix}
0_{s_h,t_1-s_h} & W_h
\end{bmatrix} \in \C^{s_h \times t_1}, \qquad W_h=\begin{bmatrix}
Y_h & \star\\
0 & \star
\end{bmatrix} \] where both $W_h \in \C^{s_h \times s_h}$ and $Y_h \in \C^{(s_h-s_{h+1}) \times (s_h -s_{h+1})}$ are either invertible upper  triangular or empty (in the case of $Y_h$ when $s_{h}=s-{h+1}$); here by convention $s_k=0$ and $\star$ denotes blocks whose precise nature is irrelevant. As consequence, denoting by $\mathrm{cwRev} N(\la)$ the columnwise reversal \cite[Section 4]{NP} of $N(\la)$, 
\[ \mathrm{cwRev} N(0) = \begin{bmatrix}
I_{t_2-t_1}\\
0
\end{bmatrix} \oplus \begin{bmatrix}
Y_1\\
0
\end{bmatrix} \oplus \cdots \oplus \begin{bmatrix}
Y_{k-2}\\
0
\end{bmatrix} \oplus W_{k-1}.\]
 Manifestly $\mathrm{cwRev}N(0)$ has full column rank, and hence $N(\la)$ is column reduced.
\end{proof}

\medskip


We show here also how to solve for the right inverse of
the same pencil, since it involves the same block Toeplitz matrix. It gives a practical flavour to the theoretical results that we obtained in Proposition \ref{prop:extractrootpoly}. We now have to solve 
\begin{equation} \label{rightinv}
\left[ \begin{array}{cccc}
E_\epsilon &  & \\
A_\epsilon & \ddots &  \\
& \ddots & E_\epsilon \\
& & A_\epsilon
\end{array}\right] \left[ \begin{array}{c}
R_{k-1}\\
\vdots \\
R_1 \\
R_0
\end{array}\right] 
=\left[ \begin{array}{c}
0 \\
\vdots \\
0 \\
I_{n}
\end{array}\right].
\end{equation}
The generalized inverse of the block diagonal matrix $A_\epsilon$ yields the minimum norm solution of the first matrix
$$  R_0= A_\epsilon^{\dagger} = \left[\begin{array}{cccc}
 0 & \ldots & 0 \\ A_{1,2}^{-1} &  \ddots & \vdots \\
 &  \ddots  & 0 \\
& & A_{k-1,k}^{-1} 
\end{array} \right].
$$
and the recurrence $A_\epsilon R_{i-1}+E_\epsilon R_i=0$ then yields the next matrices $R_{i-1}=-A_\epsilon^{\dagger} E_\epsilon R_i$.
We can again use the submatrices $Z_i$ given in \ref{Zi} to find an expression for $Z:=-A_\epsilon^{\dagger}E_{\epsilon}$
$$  Z:= -A_\epsilon^{\dagger}E_{\epsilon} = \left[\begin{array}{cccc}
 0 & \ldots & \ldots & 0 \\ Z_1 &  \ddots & & \vdots \\
 &  \ddots  & \ddots & 0 \\
& & Z_{k-1} & 0
\end{array} \right],
$$
which is a nilpotent matrix of degree $k-1$. Therefore the recurrence stops with $R_k=0$.

\begin{corollary}
Let the pencil $A_\epsilon + \la E_\epsilon$ be as in Theorem \ref{th:rightinv}. Then 
its right inverse is given by the polynomial matrix 
$$  A_\epsilon^{\dagger}\left[\sum_{i=0}^k Z^i(-\la)^i\right],
$$
where 
$$  A_\epsilon^{\dagger} = \left[\begin{array}{cccc}
 0 & \ldots & 0 \\ A_{1,2}^{-1} &  \ddots & \vdots \\
 &  \ddots  & 0 \\
& & A_{k-1,k}^{-1} 
\end{array} \right], \quad  Z:= \left[\begin{array}{cccc}
 0 & \ldots & \ldots & 0 \\ Z_1 &  \ddots & & \vdots \\
 &  \ddots  & \ddots & 0 \\
& & Z_{k-1} & 0
\end{array} \right], \quad  Z_i=\left[\begin{array}{cc} 0 & - A_{i,i+1}^{-1}\hat E_{i,i}\end{array}\right].
$$
\end{corollary}

Note that solving for the right inverse directly using the block Toeplitz equation
\eqref{rightinv} would avoid the cumbersome block bidiagonalization. However, the calculation of the minimum norm solution would then be of higher complexity.
 
\subsection{Constructing the root polynomials}

In this section we look at the calculation of a maximal set of root polynomials for the pencil $A_{0}+\la E_{0}$ with only elementary divisors at 
$\la =0$. By Proposition \ref{prop:extractrootpoly} and Theorem \ref{thm:extractrootpolys}, this suffices for the computation of a maximal set for the whole pencil \eqref{eq:goodform}. We suppose that we have  constructed the following bidiagonal staircase form using the techniques described in Section \ref{Sec:Extraction} 
\begin{equation}\label{eq:pencil}
A_0 + \la E_0=
\left[\begin{array}{cccc}
 \la  E_{1,1} & A_{1,2} &   &  \\
 &  \ddots &  \ddots &  \\
 &  & \la E_{k-1,k-1}  & A_{k-1,k} \\
 &  &  &  \la E_{k,k}
\end{array} \right] ,
\end{equation}
where 
$$ A_{i,i+1} = \left[\begin{array}{c}\hat A_{i,i+1} \\ 0 \end{array}\right] \in \C^{s_i\times t_{i+1}}, \quad E_{i,i}\in \C^{s_i\times t_i}, 
$$
and both $\hat A_{i,i+1}$ and $E_{i,i}$ are invertible upper triangular matrices. Observe that this implies that $t_{i}=s_i$ for $i=1:k$. Moreover, as noted in previous sections, there are precisely $t_i - t_{i+1}$ (recall that by convention $t_{k+1}=0$) partial multiplicities equal to $i$ for all $i=1,\dots,k$.

As pointed out earlier, the condition $(A_0+\la E_0)\sum_{i=0}^{k-1} x_i\la^i \equiv 0 \mod \lambda^k$ is 
equivalent to the convolution equation
$$ \left[ \begin{array}{cccc}
A_0 & E_0 &  & \\
& A_0 & \ddots &  \\
& & \ddots & E_0 \\
& &  & A_0
\end{array}\right] \left[ \begin{array}{c}
x_{k-1}\\
\vdots \\
x_1 \\
x_0
\end{array}\right] =
0.
$$
In order to have all such solutions with $x_0\neq 0$ we look for a block version of the solution 
$\sum_{i=0}^{k-1} X_i\la^i$ with $X_0$ of full column rank. But rather than finding a particular 
nullspace of a large block Toeplitz matrix, we exploit the bidiagonal staircase form to construct the solutions. 
Let $\hat Z_i$ be the solution of
$$  E_{i,i}\hat Z_i + A_{i,i+1} =0, \quad \Rightarrow \; \hat Z_i=-E_{i,i}^{-1}A_{i,i+1}=-E_{i,i}^{-1} \left[\begin{array}{cc}  \hat A_{i,i+1} \\ 0 \end{array}\right] \in\C^{s_{i}\times t_i}.
$$
It is then easy to see that the each column of the $n\times t_i$ block vectors below is a root polynomial of order $i$ for all $i=k,k-1,\dots,2,1$:
$$  \left[\begin{array}{c} \hat Z_1 \cdots \hat Z_{k-1}\\ \la \hat{Z}_2 \cdots \hat{Z}_{k-1}\\ \vdots \\ \la^{k-2} \hat Z_{k-1} \\ \la^{k-1} I_{t_k} \end{array}\right], \quad \dots \quad \left[\begin{array}{c} \hat Z_1 \\ \la I_{t_2}\\0\\ \vdots \\ 0 \end{array}\right], \quad \left[\begin{array}{c} I_{t_1}\\ 0 \\ \vdots \\ \vdots \\ 0 \end{array}\right] .
$$
For all triples of integers $1 \leq b, c \leq a$, let us now introduce the notation $I_{a,b:c} = \begin{bmatrix}
e_b & \dots & e_c
\end{bmatrix} \in \C^{a \times (c-b+1)}$, with the convention that such a matrix is empty if $c < b$. We now extract the rightmost $t_i - t_{i+1}$ columns from each of the block vectors above, yielding
\begin{equation}\label{eq:rp}
\begin{bmatrix}
\hat{Z}_1 \dots \hat{Z}_{k-1}\\
\la \hat{Z}_2 \cdots \hat{Z}_{k-1}\\
\vdots\\
\la^{k-2} \hat{Z}_{k-1}\\
\la^{k-1}I_{t_k}
\end{bmatrix} \in \C[\la]^{n\times t_k}, \quad \dots \begin{bmatrix}
\hat{Z}_1 I_{t_2,t_3+1:t_2}\\
\la I_{t_2,t_3+1:t_2}\\
0\\
\vdots\\
0
\end{bmatrix} \in \C[\la]^{n\times (t_2-t_3)} , \quad  
\begin{bmatrix}
I_{t_1,t_2+1:t_1}\\
0\\
\vdots\\
\vdots\\
0
\end{bmatrix} \in \C[\la]^{n\times (t_1-t_2)}.
\end{equation}

\begin{theorem} Let $A_0 + \la E_0$ be a bidiagonal pencil
$$ A_0 + \la E_0=
\left[\begin{array}{cccc}
 \la  E_{1,1} & A_{1,2} &   &  \\
 &  \ddots &  \ddots &  \\
 &  & \la E_{k-1,k-1}  & A_{k-1,k} \\
 &  &  &  \la E_{k,k}
\end{array} \right] ,
$$
where 
$$ A_{i,i+1} = \left[\begin{array}{c}\hat A_{i,i+1} \\ 0 \end{array}\right] \in \C^{s_i\times t_{i+1}}, \quad E_{i,i}\in \C^{s_i\times t_i}, 
$$
and both $\hat A_{i,i+1}$ and $E_{i,i}$ are invertible upper triangular matrices.
The its only root is $\la=0$ and the columns of the block vectors in \eqref{eq:rp} are a maximal set of root polynomials at $0$.
\end{theorem}

\begin{proof}
It is clear that each column is a root polynomial, since they are a subset of a larger set of root polynomials. We now proceed by steps:
\begin{enumerate}
\item  The columns of the block vectors in \eqref{eq:rp} are a $0$-independent set since if we put them next to each other to form a polynomial matrix and evaluate it at $0$ we obtain
\[ \begin{bmatrix}
\hat{Z}_1 \cdots \hat{Z}_{k-1} & \dots & \hat{Z}_1 I_{t_2,t_3+1:t_2}  & I_{t_1,t_2+1:t_1} \\
0 & 0 & \dots & 0\\
\vdots & \ddots & \ddots & \vdots \\
0 & \dots & 0 & 0
\end{bmatrix} \in \C^{n \times t_1};\]
moreover, the top $t_1 \times t_1$ block of the latter matrix is, by construction, invertible upper triangular, and hence its columns are linearly independent. (Note that \eqref{eq:pencil} is regular.)
\item The columns of the block vectors in \eqref{eq:rp} are a complete set, since there are precisely $t_1$ nonzero partial multiplicities of $0$ in the pencil \eqref{eq:pencil}.
\item  Finally, the columns of the block vectors in \eqref{eq:rp} are a maximal set by Theorem \ref{thm:frovan}: indeed, their orders correspond precisely to the partial multiplicities of $0$ as an eigenvalue of \eqref{eq:pencil}.
\end{enumerate} 
\end{proof}

\section{Numerical aspects} \label{Sec:Numerics}

\subsection{A worked out example}
In this subsection we illustrate our procedure using a staircase form with the zero and non-zero pattern of Example \ref{matlabex}. We generated ten such random pencils $\tilde A+\la \tilde E$ using the Matlab function
{\tt randn}. We then normalized the pencil such that $\max(\|\tilde A\|_2,\|\tilde E\|_2)=1$ and ran the bidiagonalization algorithm described in Section \ref{Sec:Extraction}.
All computations were performed with Matlab R2020a on a laptop with machine epsilon $\epsilon \approx 2\cdot 10^{-16}$. Rather than computing the nullspace and root polynomials via the blocks $Z_i$ of the bidiagonal form, we  (equivalently) reduced the bidiagonal form further to the  permuted Kronecker form corresponding to the staircase form of the pencil. For 
the Example \ref{matlabex}, this would be
\begin{equation} \label{Kron}
A_K + \la E_K := \left[\begin{array}{ccccc|ccc|c}
0 & \la & 0 & 0 & 0 & 1 & 0 & 0 & 0 \\
0 & 0 & \la & 0 & 0 & 0 & 1 & 0 & 0 \\
0 & 0 & 0 & \la & 0 & 0 & 0 & 1 & 0 \\
0 & 0 & 0 & 0 & \la & 0 & 0 & 0 & 0 \\ \hline
0 & 0 & 0 & 0 & 0 & 0 & \la & 0 & 1 \\
0 & 0 & 0 & 0 & 0 & 0 & 0 & \la & 0 
\end{array}\right]:=
S^{-1}(\tilde A+\la \tilde E)T.
\end{equation}
This is still a bidiagonal form but now with non-zero triangular blocks that have been transformed to identity
matrices of matching dimensions. 
Moreover, this Kronecker-like form can be obtained by applying non-singular upper triangular transformation matrices $S$ and $T$.  This pencil is in its Kronecker canonical form, up to a row and column permutation, and the calculation of its root polynomials $r_i(\la)$ and nullspace vectors $n_i(\la)$ are then trivial. 
They are given by 
$$  r_1=\left[\begin{array}{c}0\\0\\0\\0\\1\\ \hline 0\\0\\0\\ \hline 0 \end{array}\right] , \;
 r_2=\left[\begin{array}{c}0\\0\\0\\1\\0\\ \hline 0\\0\\ -\la \\ \hline 0 \end{array}\right], \;
 n_1=\left[\begin{array}{c}1\\0\\0\\0\\0\\ \hline 0\\0\\0\\ \hline 0 \end{array}\right] , \;
 n_2=\left[\begin{array}{c}0\\1\\0\\0\\0\\ \hline -\la \\0\\ 0 \\ \hline 0 \end{array}\right], \;
 n_3=\left[\begin{array}{c}0\\0\\1\\0\\0\\ \hline 0\\0\\0\\ \hline \la^2 \end{array}\right].
$$
It follows from \eqref{Kron} that the corresponding vectors $\tilde r_i(\la)$ and $\tilde n_i(\la)$ of the pencil $\tilde A + \la \tilde E$ are then just given by $\tilde r_i(\la)=Tr_i(\la)$ and $\tilde n_i(\la)=Tn_i(\la)$.

{\small
\begin{table}[h] 
\begin{center}
\caption{Calculation of the root vectors and null space vectors}
\begin{tabular}{cccccccc} \label{T1}
$\epsilon\cdot \kappa$ & Back & Off & $\|ResN\|$ & $\|N\|$ & $\|ResR\|$ & $\| R\|$ \\ \hline
   2.2651e-14  & 2.8868e-15  & 1.2617e-10  & 1.6326e-14  & 2.1867e+04 &  5.6850e-14 &  8.5175e+05 \\
   2.4177e-14  & 6.1009e-16  & 7.6841e-12  & 3.4516e-15  & 7.9996e+03 &  5.3857e-15 &  4.0791e+04 \\
   2.6777e-14  & 2.0687e-14  & 5.0494e-12  & 1.1192e-15  & 1.1260e+03 &  3.1776e-14 &  7.0917e+05 \\
   3.7079e-14  & 3.8283e-14  & 1.0437e-10  & 4.2609e-15  & 2.5840e+03 &  1.7053e-13 &  1.8444e+05 \\
   7.1896e-14  & 6.8027e-16  & 2.9682e-11  & 5.8754e-15  & 9.9160e+02 &  7.2462e-14 &  6.9345e+04 \\
   7.6699e-14  & 1.1937e-14  & 3.2024e-11  & 1.1879e-16  & 5.9527e+02 &  1.7798e-15 &  8.7368e+04 \\
   8.0162e-14  & 6.3404e-15  & 2.9494e-09  & 9.0764e-17  & 5.0684e+01 &  2.2204e-16 &  7.0384e+03 \\
   1.0326e-13  & 5.7180e-15  & 5.3124e-12  & 6.6572e-16  & 3.5718e+03 &  1.4211e-14 &  4.3890e+05 \\
   3.6364e-13  & 1.8848e-15  & 4.3257e-11  & 1.0596e-15  & 1.0001e+02 &  2.4882e-14 &  2.4401e+04 \\
   4.5965e-13  & 2.1778e-16  & 1.8087e-11  & 1.5271e-15  & 8.1174e+03 &  1.9214e-14 &  1.4896e+06
\end{tabular}  
\end{center}
\end{table}
}
In Table \ref{T1}, we tabulate the following quantities. Letting $\epsilon$ be the machine precision, the first column gives $\epsilon \cdot \kappa := \epsilon \cdot (\|S\|_2\|T^{-1}\|_2)$, which indicates what we can expect as error level when applying the equivalence transformation  $S(\hat A_K+\la \hat E_K)T^{-1}$ on the computed pencil 
$\hat A + \la \hat E_K$. In order to estimate the backward errors of the bidiagonalization, we show the norm of the residual pair $$ (S \hat A_K T^{-1}-\tilde A , S \hat E_K T^{-1} - \tilde E).$$
This can be considered as the backward error of the bidiagonalization step and its norm is denoted by Back. It follows from the given data that this backward error is very reasonable. The so-called off-norm Off is defined as the sum of the norms of the pair 
$(A_K-\hat A_K, E_K-\hat E_K)$. It can be seen that the Off norm of the reduction to bidiagonal form is not of the order of the machine precision. Therefore, iterative refinement should probably be applied.
 The matrices $\tilde N$ and $\tilde R$ have the null vectors $\tilde n_i(\la)$ and root vectors $\tilde r_i(\la)$ as columns
and their norm is the Frobenius norm. The residual matrices $ResN$ and $ResR$ have the vectors 
$(\tilde A+\la \tilde E)\tilde n_i(\la)$ and $(\tilde A+\la \tilde E)\tilde r_i(\la)\mod \la^k_i$ as columns. Both these matrices
are zero when the null space and root vectors are computed exactly. Our experiments show that the residuals are very close to the machine precision, despite the fact that the norms of $N$ and $R$ are quite large.

\subsection{A note on the complexity of the algorithm} 

The complexity of the method proposed in this paper for computing the null space and root polynomials of an $m\times n$ matrix pencil, is cubic in the dimensions of the pencil, i.e.\ it is of the order of ${\mathcal O}(m+n)^3$, when using an appropriate implementation. 
The most time consuming step is the orthogonal reduction to staircase form, and this was shown to have cubic complexity in \cite{BeeV88},
provided one uses orthonormal transformations to echelon form for the basic steps. The reduction to bidiagonal form requires the construction, and multiplication with triangular matrices $S$ and $T$
of respective dimensions $m$ and $n$, which again has cubic complexity. A step of iterative refinement repeats the same triangular elimination, but with a different right hand side
and is therefore also of cubic complexity. The matrices $Z_i$ and their nested products, are implicitly computed when further reducing the pencil to Kronecker-like form, which again has cubic complexity.
The vector coefficients of the polynomial vectors of the nullspace
and of the root polynomials are then vectors to be extracted from the transformation matrices $S$ and $T$, and we never have to manipulate vectors of growing dimensions, such as in the methods that use block Toeplitz matrix equations of growing dimensions. The methods that
do exploit these block Toeplitz equations typically have a complexity of the order of ${\mathcal O}(m+n)^3d^3$ where $d$ is the degree of the largest nullspace vector.



%

\section{Conclusions} \label{Sec:Conclusion}

In this paper we have devised a numerical method to compute both a maximal set of root polynomials  at $\la_0$
and a minimal basis for the null space of a given pencil. The method is based on three basic steps:
\begin{enumerate}
\item first we apply a unitary equivalence transformation to put the pencil in a particular staircase form
 which displays the right minimal indices of the pencil, and the Segr\'e characteristic of the eigenvalue
$\la_0$.
\item then we perform a block upper triangular equivalence transformation (and if necessary also a permutation) that yields a block bidiagonal
pencil and moreover separates the Kronecker part from Jordan structure at $\la_0$, and
\item next we construct via simple recurrences  or via a further reduction to Kronecker-like form the requested root polynomials and minimal basis vectors
\item finally we put everything together to obtain both a minimal basis and a maximal set of root polynomials for the original pencil.
 \end{enumerate}
At least some of these individual steps are known to be numerically stable; a rigorous analysis of the stability of other steps, and of the algorithm as a whole, is an interesting problem for future research. Namely:

\begin{enumerate}
\item In the first step we only perform unitary transformations and the backward stability of this operation has been shown in \cite{vd79}.
\item The second step can be interpreted as 
just back substitution for the solution of a set of linear equations and by using iterative refinement on this system.
\item The third step involves the solutions $Z_i$ of particular systems of equations, 
and the calculation of their products
 or the  further reduction to Kronecker-like form. 
\item The last step is trivial to implement in a stable manner as we have shown that we just need form block vectors where some of the blocks are the previously computed minimal basis or root polynomials, and every other block is zero.
\end{enumerate}   
It may be possible to implement each step at least in a forward stable manner; or, if necessary, one may even run the iterative refinement step in extended precision to guarantee a sufficiently small error. Nevertheless, the composition of (forward) stable algorithms is stable only under certain conditions \cite{bnv}. Thus, a more in-depth analysis of the stability of the whole process is of course a more subtle issue and a potential subject for future research.

\section*{Acknowledgements} We thank two anonymous reviewers, whose very insightful comments improved the paper.

\end{document}